\documentclass[12pt,reqno]{article}
\usepackage{amsmath,amsfonts,amssymb,amsthm,indentfirst}

\catcode`\@=11

 \def\AMSTeXfeatures{\Plainheads 
   \let\current@vert=\AMS@vert}

 \def\Plainheads{\sh@ftdiam=0.05em
   \getlabeldims
   \let\vshaftfill=\plnvsolidfill
   \let\hshaftfill=\plnhsolidfill
   \let\th@rhead=\plnrhead
   \let\th@lhead=\plnlhead
   \let\th@dnhead=\plndnhead
   \let\th@uphead=\plnuphead}
 
 \def\glet{\global\let}

 \def\LaTeXfeatures{\catcode`\@=11
   \ifx\@clnwd\undefined \nol@g
      \input ltxcode.tex \dol@g \fi
   \ltxheads \let\current@vert=\new@vert
   \providelto \catcode`\@=\active}

 \def\nol@g{\def\wlog{\edef\garbage}}
 \def\dol@g{\let\wlog=\wl@g} \let\wl@g=\wlog
 \nol@g % Silences allocations with \newbox, etc.

 \newbox\ltobox
 \def\providelto{{\setbox\z@=
   \hbox{$\to$}\minharrlen=\wd\z@
   \global\setbox\ltobox=\hbox{$\activeat>>>$}}
   \def\lto{\mathrel{\copy\ltobox}}}

 \def\ltxheads{\sh@ftdiam=\@wholewidth
   \getlabeldims
   \let\vshaftfill= \ltxvsolidfill
   \let\hshaftfill=\ltxhsolidfill
   \let\th@rhead=\ltxrhead
   \let\th@lhead=\ltxlhead
   \let\th@dnhead=\ltxdnhead
   \let\th@uphead=\ltxuphead}
 {\catcode`\@=\active
   \gdef@#1{\csname #1\string@at\endcsname}
   \glet\activeat=@}
 \def\def@#1{\expandafter\def\csname #1@at\endcsname}

 \def@>#1>#2>{\@rrow R{#1}{#2}}
 \def@<#1<#2<{\@rrow L{#1}{#2}}
 \def@ V#1V#2V{\@rrow V{#1}{#2}}
 \def@ A#1A#2A{\@rrow A{#1}{#2}}
 \def@/#1/#2/#3/{\@rrow{#1}{#2}{#3}}
 \def@.{\ifodd\row\ifmmode\noharrow
     \else\leavevmode.\spacefactor3000 \fi
   \else\novarrow\fi}
 \def@={\ifodd\row\harrow\hequalfill{}{}%
   \else\varrow\vequalfill{}{}\fi}
 \def@:#1{\ifx=#1\harrow\deffill{}{}%
   \else\leavevmode\null:#1\fi}
 \def@|{\current@vert}
  \def\AMS@vert{\varrow\vequalfill{}{}}
  \def\new@vert#1|#2|{\ifodd\row
   \let\nextarrow\vertexvarrow
   \else\let\nextarrow\varrow\fi
   \nextarrow\vshaftfill{#1}{#2}}
 \def@-{\ifmmode\let\next\hl@ne
   \else\let\next\AMSatdash \fi \next}
  \def\hl@ne#1-#2-{\harrow\hshaftfill{#1}{#2}}
  \def\AMSatdash{\let\next\relax\leavevmode
    \def\next@{\ifx\next-%
      \def\next-{\futurelet\next\nextii@}%
     \else\def\next{\hbox{-}}\fi\next}%
    \def\nextii@{\ifx\next-\def\next-{\hbox{---}}%
      \else\def\next{\hbox{--}}\fi\next}%
    \futurelet\next\next@}
 \def@(#1){\tweenarrows{#1}}
 \def@[#1]{\setsp@n#1\relax\activeat}
 \def\fiberbox{\hbox{$\vcenter{\hr@le\hbox{\vr@le
   \kern1ex\vbox{\kern1.2ex}\vr@le}\hr@le}$}}
  \def\hr@le{\hrule height \sh@ftdiam}
  \def\vr@le{\vrule width \sh@ftdiam}
 \def@+#1+#2+#3+{\ifodd\row \harrow{#1}{#2}{#3}%
   \else \varrow{#1}{#2}{#3}\fi}

 \def\Dnarrfill{\vequalfill\Dnhe@d}
 \def\Uparrfill{\Uphe@d\vequalfill}
 \def\ontofill{\rtarrfill\kern-0.3em %2\he@dwd
   \th@rhead\kern 0.3em} %new def

 \def\rtarrfill{\hshaftfill\th@rhead}
 \def\ltarrfill{\th@lhead\hshaftfill}
 \def\dnarrfill{\vshaftfill\th@dnhead}
 \def\uparrfill{\th@uphead\vshaftfill}
 \def\hequalfill{\plnhfill=}
 \def\deffill{:\plnhfill=}
 \def\plnvextfill#1{\setbox\z@
   \hbox{\the\textfont3 #1}%
   \dimen@=\dp\z@\advance\dimen@\ht\z@
   \copy\z@ \kern-\dimen@ %-\dp\z@
   \cleaders\copy\z@ \vfill
   \kern-\dimen@ %-\dp\z@
   \box\z@}
 \def\plnhfill#1{$\m@th\mkern-1.5mu\mathord#1\mkern-6mu
    \cleaders\hbox{$\mkern-2mu\mathord#1\mkern-2mu$}\hfill
    \mkern-6mu\mathord#1\mkern-1.5mu$}
 \def\vequalfill{\plnvextfill{\char'167}}
 \def\plnvsolidfill{\plnvextfill{\char'077}}
 \def\plnhsolidfill{\plnhfill-}
 \def\ltxhsolidfill{\leaders\hrule height\topofshaft depth\botofshaft
   \hfill}
 \def\ltxvsolidfill{\leaders\vrule width\sh@ftdiam\vfill}
 \def\hdashfill{\hd@sh\wd@sh
   \xleaders \hbox{\wd@sh\hd@sh\wd@sh}\hfill
   \wd@sh\hd@sh}
 \def\vdashfill{\vd@sh\wd@sh
   \xleaders \vbox{\wd@sh\vd@sh\wd@sh}\vfill
   \wd@sh\vd@sh}
 \def\dashed{\ifinmeasureCD\else
    \ifodd\row\option{\let\hshaftfill=\hdashfill}%
   \else\option{\let\vshaftfill=\vdashfill}\fi\fi}
 \newdimen\CDstrutht  \newdimen\CDstrutdp
 \newdimen\CDstrutlen \CDstrutlen=\CDstrutht
   \advance\CDstrutlen by \CDstrutdp
 \def\CDstrut{\vrule
   height \ifnum\row=1 \z@\else\CDstrutht \fi
   depth \ifnum\row=\numrows \z@ \else\CDstrutdp \fi
   width\z@}
 \newdimen\CDarrsurr \CDarrsurr=0.375em
 \newdimen\CDdashlen
 \newdimen\CDvarrlen \CDvarrlen=1.5\baselineskip
 \newdimen\minharrlen %Used outside CD's
  \setbox\z@\hbox{$\longrightarrow$} \minharrlen=\wd\z@
 \newdimen\minCDharrlen \minCDharrlen=2.5em %825 2pc %2.5pc
\newdimen \minc@lwd
\def\findminc@lwd{\minc@lwd=2\CDarrsurr
  \advance\minc@lwd\minCDharrlen}
 \newdimen\sh@ftdiam

 \newdimen\labelsurr \labelsurr=1.25 em

\newcount\sp@ncnt \sp@ncnt=\@ne
\newcount\sp@ncnt@ \sp@ncnt@=\@ne
\newdimen\@rrwd \newdimen\@rrdp

 \def\adjustbot#1{\option{\advance\@rrdp#1\relax}}
\def\pushvertex#1{\global\p@shlen#1\relax
   \global\let\maybepush=\dopush}

 \newdimen\p@shlen \p@shlen=\z@

 \let\maybepush=\relax
 \def\dopush{\ifinmeasureCD %omitted by accident
   \advance\locdimen by -\p@shlen %AL
   \else\advance \@rrwd by -\p@shlen \fi %AL
   \global\let\maybepush=\relax \global\p@shlen=\z@\relax}
 \def\span@ne{\global\sp@ncnt=\@ne\relax}
 \def\setsp@n#1#2{\global\sp@ncnt=#1\relax
   \ifx\relax#2\relax\else\global\sp@ncnt@=#2\relax\fi}

 \def\plnrhead{\llap{$\rightarrow\mkern-1.5mu$}}
 \def\plnlhead{\rlap{$\mkern-1.5mu\leftarrow$}}

 \def\clap#1{\hbox to \z@{\hss #1\hss}}

 \def\plndnhead{\hbox{\the\textfont3 \char'171}}
 \def\plnuphead{\hbox{\the\textfont3 \char'170}}
 \def\Dnhe@d{\hbox{\the\textfont3 \char'177}}
 \def\Uphe@d{\hbox{\the\textfont3 \char'176}}

 \def\ltxrhead{\raise\@xisheight
   \llap{\smash{\@linefnt\@getrarrow(1,0)}}}
 \def\ltxlhead{\raise\@xisheight
   \rlap{\@linefnt\@getlarrow(-1,0)}}
 \def\ltxuphead{\setbox\z@=\rlap{%
   \kern\@halfwidth\@linefnt\char'66}%
   \copy\z@\kern-\ht\z@}
 \def\ltxdnhead{\setbox\z@=\rlap{%
   \kern\@halfwidth\@linefnt\char'77}%
   \ht\z@=\z@\box\z@}

 \def\wd@sh{\kern0.5\CDdashlen}
 \def\hd@sh{\vrule height\topofshaft depth\botofshaft
    width\CDdashlen}
 \def\vd@sh{\hrule height\CDdashlen
   depth\z@ width\sh@ftdiam}
\def\xylist{14{3434}13{2414}12{1723}%
  23{1413}34{1153}11{0867}43{0707}%
  32{0580}21{0414}31{0291}41{0}}
\newcount\tgtcnt@
\def\find@xyargs{\dimen@=\@rrdp
  \advance\dimen@ by \CDstrutlen
  \tgtcnt@=\dimen@ \dimen@=\@rrwd %\relax
  \divide\dimen@ by \@m %\relax
  \divide \tgtcnt@ by \dimen@ %\relax
  \expandafter\testxy\xylist\relax
  \unitlength=\@xarg\@rrdp
  \divide\unitlength by\@yarg\relax}
\def\testxy#1#2#3{\ifnum\tgtcnt@>#3
    \@xarg=#1\relax \@yarg=#2\relax
    \let\next=\ignorerest
  \else\let\next\testxy\fi\next}
\def\ignorerest#1\relax{\relax}

\let\scalefactor=\@ne
\def\SWarrow{\find@xyargs\vector
  (-\@xarg,-\@yarg)\scalefactor\hskip-\wd\@linechar}
\def\NWarrow{\find@xyargs\vector
  (-\@xarg,\@yarg)\scalefactor\hskip-\wd\@linechar}
\def\NEarrow{\find@xyargs\vector
  (\@xarg,\@yarg)\scalefactor}
\def\SEarrow{\find@xyargs\vector
  (\@xarg,-\@yarg)\scalefactor}
\def\rightupline{\find@xyargs\@linelen=\scalefactor
     \unitlength\@sline}
\def\rightdownline{\find@xyargs\@yarg=-\@yarg\relax
     \@linelen=\scalefactor\unitlength\@sline}

\def\Sim{\ifodd\row\setbox\z@=\hbox{$\sim$}\dimen@=\ht\z@
 \advance\dimen@ by -\@xisheight
  \vbox{\box\z@\kern-\@xisheight\kern\dimen@}%
  \else\hbox{$\wr$}\fi}

\def\harrow#1#2#3{\inmeasureCDtrue\findminarrwd
  {#2}{#3}{\sp@ncnt\minharrlen}\inmeasureCDfalse\span@ne
  \mathrel{\hbox{\options\hplace{#1}\ulabel{#2}\dlabel{#3}}}}

\def\noharrow{\harrow\hfill{}{}}
\def\vertexvarrow#1#2#3{\findarrdp \@rrwd=\z@ \setsp@n\@ne\@ne
  \vbox to \z@{\kern-1.2\CDstrutht
  \rlap{\options\vplace{#1}\llabel{#2}\rlabel{#3}}\vss}}

\newif\ifinmeasureCD
\def\measurelabel#1{\setbox\z@
  \hbox{$\scriptstyle#1\kern\labelsurr$}%
  \ifdim\wd\z@>\@rrwd \@rrwd=\wd\z@\fi}
\def\findminarrwd#1#2#3{\@rrwd=#3\relax
   \measurelabel{#1}\measurelabel{#2}}
\def\findCDarrwd#1#2{\@rrwd=\minCDharrlen
   \measurelabel{#1}\measurelabel{#2}%
  }

\newcount\row \row=\@ne \newcount\col \col=\@ne %96 initialized
 \newcount\numrows 
\numrows=\@ne
 \newcount\numcols
\newcount\arrspan \newdimen\vrtxhalfwd  \newbox\tempbox

\def\DANABUG{\advance\col by \@ne
 \@rrwd=\minCDharrlen
  \advance\@rrwd by \vrtxhalfwd
  \advance\@rrwd by \CDarrsurr
  \ifnum\col>\numcols \numcols=\col
     \newlocdimen{col\the\col}\locdimen=\@rrwd %AL
  \else \ifdim\@rrwd>\c@l \c@l=\@rrwd\fi\fi}

\def\drop#1\\{%\noharrow %caused by DANABUG
  \findvrtxhalfsum\DANABUG\advance\row by 2 \measureinit}

\def\measureinit{\col=\@ne \vrtxhalfwd=-\CDarrsurr\arrspan=\@ne\@rrwd=\z@
   \setbox\tempbox=\hbox\bgroup$}
\def\measure{%CR \bgroup
  \let\harrow\measureCDarrow
  \let\CDCR=\measureCR %CR
   \findminc@lwd 
  \inmeasureCDtrue
  \row=\@ne \numcols=\z@ \measureinit}

\def\endmeasure{\findvrtxhalfsum\DANABUG
  \numrows=\row %CR \egroup
  \inmeasureCDfalse}

\def\newlocdimen#1{\advance\dimenc@unt by \@ne
  \ifnum\dimenc@unt<\insc@unt
     \else\errmessage{No room for the CD}\fi
  \dimendef\locdimen=\dimenc@unt
  \expandafter\dimendef\csname#1\endcsname=\dimenc@unt}

 \def\r@wc@l{\csname row\the\row col\the\col\endcsname}
 \def\c@l{\csname col\the\col\endcsname}

 \def\findvrtxhalfsum{$\egroup
  \newlocdimen{row\the\row col\the\col}%%AL
  \locdimen=\vrtxhalfwd %AL
  \vrtxhalfwd=0.5\wd\tempbox %\maybes@ve %8231
  \advance\vrtxhalfwd by \CDarrsurr
  \advance\locdimen by \vrtxhalfwd %AL
  \advance\@rrwd by \locdimen %AL
  \maybepush
  \divide\@rrwd by \arrspan\relax
  \ifdim\@rrwd<\minc@lwd
    \ifnum\col>\@ne \@rrwd=\minc@lwd\fi \fi
  \loop %94 \edef\c@l{\csname col\the\col\endcsname}
    \ifnum\col>\numcols \numcols=\col
       \newlocdimen{col\the\col}% %AL
       \locdimen=\@rrwd %AL
    \else \ifdim\@rrwd>\c@l \c@l=\@rrwd\fi \fi
   \ifnum\arrspan>\@ne
      \advance\arrspan by -1 \advance\col by \@ne
  \repeat }

 \def\measureCDarrow#1#2#3{\findvrtxhalfsum
   \arrspan=\sp@ncnt\relax\global\sp@ncnt=1\relax
   \advance\col by \@ne
   \findCDarrwd{#2}{#3}%
   \setbox\tempbox=\hbox\bgroup$}

 \newcount\dr@tn \dr@tn=\z@
 \def\locate#1:#2{\ifinmeasureCD\else
   \count@=-#1
   \multiply\count@ by 2
   \advance\count@ by #2
   \dimen@=\count@\@rrwd
   \ifnum\dr@tn=\@ne\relax \else\dimen@=-\dimen@ \fi
   \dimen@i=\@rrdp
   \ifnum\dr@tn>\z@\advance\dimen@i by \CDstrutlen \fi
   \dimen@i=\count@\dimen@i
   \count@=#2 \multiply\count@ by 2
   \divide\dimen@ by \count@
   \divide\dimen@i by \count@
   \lift\dimen@i\nudge\dimen@\fi}
\def\betweenCDrows{\advance\row by \@ne \col=\@ne
\options}

\def\hbegin{\hbox\bgroup\kern\c@l \kern-\r@wc@l$}
\def\hend{$\glet\maybepush\relax \CDstrut\egroup}
\def\vbegin{\setbox\tempbox=\hbox\bgroup$}
\def\vend{$\egroup\ht\tempbox=\z@\dp\tempbox\CDvarrlen
  \box\tempbox}
\def\setCD{\let\harrow=\setCDarrow
  \let\CDCR=\setCR %CR
  \row=\@ne \col=\@ne \hbegin}
\let\endsetCD=\hend %AT (Assume CD ends with hmaterial)

\def\findarrwd{\@rrwd=\z@ \count@=\col \advance\count@ by\sp@ncnt
  \loop\ifnum\count@>\col \advance\count@ by -1
      \advance\@rrwd by\csname col\the\count@\endcsname\repeat}
\def\setCDarrow#1#2#3{\kern\CDarrsurr\advance\col by \@ne
  \findarrwd \advance\@rrwd by -\r@wc@l  
  \@rrdp=\z@ %&0211 (It might be used by \locate).
  \maybepush
  \advance\col by -\@ne \advance\col by \sp@ncnt \span@ne
  \hbox to \@rrwd{\options
   \@rrwd=\scalefactor\@rrwd\hss
   \hplace{#1}\ulabel{#2}\dlabel{#3}\hss}%
   \kern\CDarrsurr}

\newdimen\labspacei %96 use subscript min of TeXbook 13a, p.444
\newdimen\labspaceii %96 Note many letters stick down below their boxes.

\newdimen\@xisheight
  \@xisheight=\the\fontdimen22\textfont2
\newdimen\labelskip
  \labelskip=\the\fontdimen10\textfont3 %2pt
\newdimen\topofshaft
\newdimen\botofshaft
\newdimen\botofulabel
\newdimen\topofdlabel
\def\getlabeldims{
  \topofshaft=0.5\sh@ftdiam
  \botofshaft=\topofshaft
  \advance\topofshaft by \@xisheight  
  \advance\botofshaft by -\@xisheight  
  \botofulabel=\topofshaft
  \advance\botofulabel by \labelskip
  \topofdlabel=\botofshaft
  \advance\topofdlabel by \labelskip}

\def\ulabel{\ifnum\row=\@ne\let\next\ulabeli
   \else\let\next\ulabellap\fi\next}
\def\ulabeli#1{\vbox{
  \clap{\kern-\@rrwd$\scriptstyle#1$}%
  \kern\botofulabel}\maybeoffset}
\def\ulabellap#1{\vbox to \z@{\vss
  \clap{\kern-\@rrwd$\scriptstyle#1$}%
  \kern\botofulabel}\maybeoffset}
\def\dlabel{\ifnum\row=\numrows\let\next\dlabeli
   \else\let\next\dlabellap\fi\next}
\def\dlabeli#1{\vtop{\kern\topofdlabel
  \clap{\kern-\@rrwd$\scriptstyle#1$}%
  }\maybeoffset}
\def\dlabellap#1{\vbox to \z@{\kern\topofdlabel
  \clap{\kern-\@rrwd$\scriptstyle#1$}%
  \vss}\maybeoffset}
\def\rlabel#1{\vbox to \z@{\vss
  \rlap{\kern\labelskip$\scriptstyle#1$}%
  \vss\kern-\@rrdp}\maybeoffset}
\def\llabel#1{\vbox to \z@{\vss
  \llap{$\scriptstyle#1$\kern\labelskip}%
  \vss\kern-\@rrdp}\maybeoffset}
\def\swlabel#1{\vtop{\kern0.5\@rrdp
  \llap{$\scriptstyle#1$\kern\labelskip\kern-0.5\@rrwd}
  }\maybeoffset}
\def\nwlabel#1{\vbox{
  \llap{$\scriptstyle#1$\kern\labelskip\kern-0.5\@rrwd}%
  \kern-0.5\@rrdp}\maybeoffset}
\def\selabel#1{\vtop{\kern0.5\@rrdp
  \rlap{\kern0.5\@rrwd\kern\labelskip$\scriptstyle#1$}%
  }\maybeoffset}
\def\nelabel#1{\vbox{
  \rlap{\kern0.5\@rrwd\kern\labelskip$\scriptstyle#1$}%
  \kern-0.5\@rrdp}\maybeoffset}
\def\cplace#1{\vbox to \z@{\vss
  \clap{$#1$\kern-\@rrwd}%
  \kern-\@rrdp\vss}\maybeoffset}
\def\hplace#1{\hbox to \@rrwd{#1}\maybeoffset}
\def\vplace#1{\clap{\vbox to \z@{#1\kern-\@rrdp}}\maybeoffset}
\newdimen\nudgeamount \nudgeamount=\z@
\newdimen\liftamount \liftamount=\z@
\let\maybeoffset\relax
\newbox\offsetbox \newdimen\lastheight
\def\dooffset{%assumes that \lastbox is a <box> set in horiz. mode
  \setbox\offsetbox=\lastbox \lastheight=\ht\offsetbox 
  \setbox\offsetbox=\vbox{\kern-\liftamount\box\offsetbox}%
  \ht\offsetbox=\lastheight
  \kern\nudgeamount\box\offsetbox\kern-\nudgeamount
  \global\nudgeamount=\z@ \global\liftamount=\z@
  \glet\maybeoffset=\relax}
\def\nudge#1{\ifinmeasureCD\else
  \global\advance\nudgeamount#1\relax
  \global\let\maybeoffset\dooffset\fi}
\def\lift#1{\ifinmeasureCD\else
  \global\advance\liftamount#1\relax
  \global\let\maybeoffset\dooffset\fi}

\def\findarrdp{\@rrdp=\CDvarrlen
  \ifnum\sp@ncnt@>1
    \advance\@rrdp by \CDstrutlen
    \multiply\@rrdp by \sp@ncnt@
    \advance\@rrdp by -\CDstrutlen \fi
 }

\def\varrow#1#2#3{\ifnum\sp@ncnt>\@ne 
     \sp@ncnt@=\sp@ncnt\relax\fi
  \findarrdp \@rrwd=\z@ %&0211 It might be used by \locate
  \kern\c@l
   \hbox to \z@{\options
   \@rrdp=\scalefactor\@rrdp
    \hss\vplace{#1}\llabel{#2}\rlabel{#3}\hss}%
  \global\advance\col by \@ne \setsp@n\@ne\@ne
  }

\def\novarrow{\varrow\vfill{}{}}

\def\tweenarrows#1{\findarrwd \findarrdp \setsp@n\@ne\@ne
  \rlap{\options\cplace{#1}}}

\def\usarrow #1#2#3{\dr@tn=\@ne
  \findarrwd \findarrdp \setsp@n\@ne\@ne 
  \rlap{\options\cplace{#1}\nwlabel{#2}\selabel{#3}}%
  \dr@tn=\z@}
\def\dsarrow #1#2#3{\dr@tn=\tw@
  \findarrwd \findarrdp \setsp@n\@ne\@ne 
  \rlap{\options\cplace{#1}\swlabel{#2}\nelabel{#3}}%
  \dr@tn=\z@}
 \def\@rrow#1{\csname #1@rrow\endcsname}
 \def\R@rrow{\harrow \rtarrfill}
 \def\L@rrow{\harrow \ltarrfill}
 \def\V@rrow{\varrow \dnarrfill}
 \def\A@rrow{\varrow \uparrfill}
 \def\SE@rrow{\dsarrow \SEarrow}
 \def\NW@rrow{\dsarrow \NWarrow}
 \def\SW@rrow{\usarrow \SWarrow}
 \def\NE@rrow{\usarrow \NEarrow}
 \def\DS@rrow{\dsarrow \dnslope}
 \def\US@rrow{\usarrow \upslope}
 \def\upslope{\find@xyargs
       \@linelen=\unitlength\@sline}
 \def\dnslope{\find@xyargs\@yarg=-\@yarg\relax
       \@linelen=\unitlength\@sline}

\newtoks\optionlist 
\optionlist={}
\let\options\relax
\def\dooptions{\the\optionlist\global\optionlist={}%
  \glet\options=\relax}
\def\option#1{\ifinmeasureCD\else
  \glet\options=\dooptions
  \global\optionlist=\expandafter{\the\optionlist\relax#1}\fi}
\def\wider#1{\ifinmeasureCD\else
   \option{\advance\@rrwd by #1}\fi}
\def\deeper#1{\ifinmeasureCD\else
   \option{\advance\@rrdp by #1}\fi}
{\def\\{\global\let\sptoken= }\\ }%now \sptoken is a spacetoken

\def\CR{\futurelet\nexttok\testCR}
\def\testCR{\ifx\nexttok\sptoken
   \let\next\eatspaceCR\else\let\next\CDCR\fi\next}
\def\eatspaceCR#1 {\CR}
\def\measureCR{\ifx\nexttok\endmeasure\let\nextCR\relax
    \else\let\nextCR\drop\fi\nextCR}
\def\setCR{\ifodd\row
  \ifx\nexttok\endsetCD\else\hend\betweenCDrows\vbegin\fi
  \else\vend\betweenCDrows\hbegin\fi}

\countdef\dimenc@unt=11
\def\CD#1\endCD{%CRAL
   \begingroup\let\\=\CR
  \m@th\offinterlineskip
   \measure#1\endmeasure\null\,\vcenter{\setCD#1\endsetCD}\,
   \endgroup
    }

\ifx\@clnwd\undefined \nol@g\else\catcode`\ =14\relax\fi
 \font\@linefnt=line10 
 \newcount\@tempcnta
 \newcount\@tempcntb
 \newdimen\@tempdima
 \newdimen\@tempdimb
 \newdimen\@wholewidth
 \newdimen\@halfwidth
   \@wholewidth\fontdimen8\@linefnt \@halfwidth .5\@wholewidth
 \newdimen\unitlength
 \newcount\@xarg
 \newcount\@yarg
 \newcount\@yyarg
 \newbox\@linechar
 \newdimen\@linelen
 \newdimen\@clnwd
 \newdimen\@clnht
 \newif\if@negarg
 
 \def\@whilenoop#1{}

 \def\@whiledim#1\do #2{\ifdim #1\relax#2\@iwhiledim{#1\relax#2}\fi}

 \def\@iwhiledim#1{\ifdim #1\let\@nextwhile=\@iwhiledim 
         \else\let\@nextwhile=\@whilenoop\fi\@nextwhile{#1}}

 \def\@sline{\ifnum\@xarg< 0 \@negargtrue \@xarg -\@xarg \@yyarg -\@yarg
   \else \@negargfalse \@yyarg \@yarg \fi
 \ifnum \@yyarg >0 \@tempcnta\@yyarg \else \@tempcnta -\@yyarg \fi
 \ifnum\@tempcnta>6 \@badlinearg\@tempcnta0 \fi
 \ifnum\@xarg>6 \@badlinearg\@xarg 1 \fi
 \setbox\@linechar\hbox{\@linefnt\@getlinechar(\@xarg,\@yyarg)}%
 \ifnum \@yarg >0 \let\@upordown\raise \@clnht\z@
    \else\let\@upordown\lower \@clnht \ht\@linechar\fi
 \@clnwd=\wd\@linechar
 \if@negarg \hskip -\wd\@linechar \def\@tempa{\hskip -2\wd\@linechar}\else
      \let\@tempa\relax \fi
 \@whiledim \@clnwd <\@linelen \do
   {\@upordown\@clnht\copy\@linechar
    \@tempa
    \advance\@clnht \ht\@linechar
    \advance\@clnwd \wd\@linechar}%
 \advance\@clnht -\ht\@linechar
 \advance\@clnwd -\wd\@linechar
 \@tempdima\@linelen\advance\@tempdima -\@clnwd
 \@tempdimb\@tempdima\advance\@tempdimb -\wd\@linechar
 \if@negarg \hskip -\@tempdimb \else \hskip \@tempdimb \fi
 \multiply\@tempdima \@m
 \@tempcnta \@tempdima \@tempdima \wd\@linechar \divide\@tempcnta \@tempdima
 \@tempdima \ht\@linechar \multiply\@tempdima \@tempcnta
 \divide\@tempdima \@m
 \advance\@clnht \@tempdima
 \ifdim \@linelen <\wd\@linechar
    \hskip \wd\@linechar
   \else\@upordown\@clnht\copy\@linechar\fi}
 
 \def\@getlinechar(#1,#2){\@tempcnta#1\relax\multiply\@tempcnta 8
 \advance\@tempcnta -9 \ifnum #2>0 \advance\@tempcnta #2\relax\else
 \advance\@tempcnta -#2\relax\advance\@tempcnta 64 \fi
 \char\@tempcnta}
 
 \def\vector(#1,#2)#3{\@xarg #1\relax \@yarg #2\relax
 \@tempcnta \ifnum\@xarg<0 -\@xarg\else\@xarg\fi
 \ifnum\@tempcnta<5\relax
 \@linelen=#3\unitlength
 \ifnum\@xarg =0 \@vvector 
   \else \ifnum\@yarg =0 \@hvector \else \@svector\fi
 \fi
 \else\@badlinearg\fi}
 
 \def\@svector{\@sline
 \@tempcnta\@yarg \ifnum\@tempcnta <0 \@tempcnta=-\@tempcnta\fi
 \ifnum\@tempcnta <5
   \hskip -\wd\@linechar
   \@upordown\@clnht \hbox{\@linefnt  \if@negarg 
   \@getlarrow(\@xarg,\@yyarg) \else \@getrarrow(\@xarg,\@yyarg) \fi}%
 \else\@badlinearg\fi}
 
 \def\@getlarrow(#1,#2){\ifnum #2 =\z@ \@tempcnta='33\else
 \@tempcnta=#1\relax\multiply\@tempcnta \sixt@@n \advance\@tempcnta
 -9 \@tempcntb=#2\relax\multiply\@tempcntb \tw@
 \ifnum \@tempcntb >0 \advance\@tempcnta \@tempcntb\relax
 \else\advance\@tempcnta -\@tempcntb\advance\@tempcnta 64
 \fi\fi\char\@tempcnta}
 
 \def\@getrarrow(#1,#2){\@tempcntb=#2\relax
 \ifnum\@tempcntb < 0 \@tempcntb=-\@tempcntb\relax\fi
 \ifcase \@tempcntb\relax \@tempcnta='55 \or 
 \ifnum #1<3 \@tempcnta=#1\relax\multiply\@tempcnta
 24 \advance\@tempcnta -6 \else \ifnum #1=3 \@tempcnta=49
 \else\@tempcnta=58 \fi\fi\or 
 \ifnum #1<3 \@tempcnta=#1\relax\multiply\@tempcnta
 24 \advance\@tempcnta -3 \else \@tempcnta=51\fi\or 
 \@tempcnta=#1\relax\multiply\@tempcnta
 \sixt@@n \advance\@tempcnta -\tw@ \else
 \@tempcnta=#1\relax\multiply\@tempcnta
 \sixt@@n \advance\@tempcnta 7 \fi\ifnum #2<0 \advance\@tempcnta 64 \fi
 \char\@tempcnta}
\catcode`\ =10

\dol@g %925
\catcode`\@=\active
\LaTeXfeatures

\newtheorem{thm}{Theorem}
\newtheorem*{prob*}{Problem}

\newtheorem{prop}{Proposition}[section]

\theoremstyle{definition}
\newtheorem{defn}{Definition}[section]
\newtheorem*{rem*}{Remark}

\def\vp{\varphi}
\def\a{\alpha}

          \def\Aut{\textrm{Aut\,}}

    \def\Hom{\textrm{Hom}}
           \def\int{\textrm{int}}

\def\Cl{\textrm{Cl}}

\def\Grp{\textrm{Grp}}
\def\Var{\textrm{Var}}
\def\Com{\textrm{Com}}
\def\End{\textrm{End}}
\def\Rep{\textrm{Rep}}
\def\Ker{\textrm{Ker}}
\def\Id{\textrm{Id}}
\numberwithin{equation}{section}
\date{}
\title{Problems in algebra inspired by universal algebraic geometry}
\author{Boris Plotkin, Hebrew University, Jerusalem}

\begin{document}

\maketitle
\begin{abstract}
Let $\Theta$ be a variety of algebras. In every $\Theta$ and every
algebra $H$ from $\Theta$ one can consider algebraic geometry in
$\Theta$ over $H$. We consider also a special categorical
invariant $K_\Theta (H)$ of this geometry. The classical algebraic
geometry deals with the variety $\Theta=Com-P$ of all associative
and commutative algebras over the ground field of constants $P$.
An algebra $H$ in this setting is an extension of the ground field
$P$.  Geometry in groups is related to varieties $\Grp$ and
$\Grp-G$, where $G$ is a group of constants. The case $\Grp -F$
where $F$ is a free group, is related to Tarski's problems devoted
to logic of a free group.

The described general insight on algebraic geometry in different
varieties of algebras inspires some new problems in algebra and
algebraic geometry. The  problems of such kind determine, to a
great extent, the content of universal algebraic geometry.

For example, a  general and natural problem is:
\medskip

{\it When do the algebras $H_1$ and $H_2$ have the same
geometry?}
\medskip

\noindent
or more specifically,
\medskip

{\it What are the conditions on algebras from a given variety
$\Theta$ which provide coincidence of their algebraic geometries?}
\medskip

\noindent
 We consider two variants of coincidence:
\medskip

1) $K_\Theta (H_1)$ and $K_\Theta (H_2)$ are isomorphic.
\medskip

2) These categories are equivalent.
\medskip

\noindent
 This problem is highly connected with the following
general algebraic problem:

Let $\Theta^0$ be the category of all free in $\Theta$ algebras
$W=W(X)$, where $X$ is finite. Consider the groups of
automorphisms $\Aut(\Theta^0)$ for different varieties $\Theta$
and also the groups of autoequivalences of $\Theta^0$. The problem
is to describe these groups for different $\Theta$.

We start with the short overview of main definitions and results
and then consider the list of unsolved problems. The results
without  references can be found in \cite{Pl4}

%In the paper the main attention is paid to conditions on algebras
%from a given variety $\Theta$ which provide coincidence of their
%algebraic geometries.
% under which two algebras from a variety of
%algebras have the same algebraic geometry.
%  The main part here
%play the notions mentioned in the title of the paper.
\end{abstract}

\section{ Definitions}

{\bf 1.1.} Fix a variety $\Theta$.  Take an algebra $H\in \Theta$
and a free in $\Theta$ algebra $W=W(X)$ with finite $X$.  The set
of homomorphisms $\Hom (W, H)$ we consider as an affine space of
points over $H$.  Points of this space are the homomorphisms
$\mu:W\to H$. If $X = \{ x_1, \cdots, x_n\}$, then we have a
bijection
\[
\alpha_X :\Hom(W, H) \to H^{(n)},
\]
defined by $\a_X(\mu) = (\mu(x_1),\cdots,\mu(x_n))$.  A point $\mu
$ is a root of the pair $(w, w')$, $w, w'\in W$, if
$w^\mu=w'^\mu$, which means also that $(w, w') \in \Ker \mu$. Here
$\Ker \mu$ is, in general, a congruence of the algebra $W$.
Simultaneously, $\mu$ is a solution of the equation $w = w'$.  We
will identify the pair $(w, w')$ and the equation $w = w'$.
%%%%%%%%%%%%%%%%%%%%%%%%%%%%%%%%%%%%%%%%%%%%%%%%%%%%%%%%%%%%%%%%%%%%%%%%%%%%%%%%%%

%\subsection{Galois correspondence}

Let now $T$ be a system of equations in $W$ and $A$ a set of
points in $\Hom(W, H)$.  We have the following Galois
correspondence
\[
\left\{%\begin{matrix}
\begin{aligned}&T'_H =\{ \mu:W\to H\ | T\subset \Ker\mu\}\\
&A'_W =\bigcap\limits_{\mu\in A} \Ker \mu
%\end{matrix}
\end{aligned} \right.
\]

\begin{defn}
The set $A$ of the form $A=T'$ for some $T$ we call a (closed)
algebraic set.  The congruence $T$ of the form $T=A'$ for some $A$
is an $H$-closed congruence.
\end{defn}

It is easy to see that the congruence $T$ is $H$-closed if and
only if $W/T \in SC(H)$, where $S$ and $C$ are the  operators of
taking subgroups and cartesian products on group classes.

One can consider the closures $A^{''} = (A')'$ and $T_H^{''} =
(T'_H)'$.
\begin{prop}
The pair $(w_0, w'_0)$ belongs to  $T^{''}_H$ if and only if the
formula ({\it  infinitary quasiidentity})
\[
\Big(\bigwedge\limits_{((w, w')\in T} (w\equiv w')\Big)
\Rightarrow w_0 \equiv w'_0
\]
holds in $H$.
\end{prop}
%\subsection{Categories}

{\bf 1.2.} We have defined the category $\Theta^0$.  Let us add to
the definition that for all objects of $\Theta^0$ the finite $X$ are
subsets of an infinite universum $X^0$. Then $\Theta^0$ is a small
category.

Define, further, the {\it category of affine spaces}
$K^0_\Theta(H)$. Objects of this category are affine spaces
\[
\Hom(W, H),\; \;  W \in 0b\; \Theta^0.
\]
The morphisms
\[
\tilde s: \Hom(W(X), H) \to \Hom (W(Y), H)
\]
of  $K^0_\Theta(H)$  are determined by homomorphisms $s:W(Y) \to
W(X)$ by the rule $\tilde s(\nu) = \nu s$ for every $\nu: W(X) \to
H$.  We have a contravariant functor
\[
\varphi  : \Theta^0\to K^0_\Theta (H).
\]

\begin{prop}
The functor $\varphi : \Theta^0\to K^0_\Theta (H)$ determines
duality of categories if and only if $Var (H) = \Theta$.
\end{prop}
%\noindent {\bf Corollary} . If $Var(H_1) = Var (H_2) = \Theta$,
%then the categories $K^0_\Theta(H_1)$ and $K^0_\Theta(H_2)$ are
%isomorphic.

Proceed now to the {\it category of  algebraic sets}
$K_\Theta(H)$. Its objects have the form $(X, A)$, where $A$ is an
algebraic set in the space $\Hom(W(X), H)$.  The morphisms $[s]:
(X, A) \to (Y, B)$ are defined by those $s: W(Y) \to W(X)$, for
which $\tilde s(\nu) \in B$ if $\nu \in A$.   Simultaneously, we
have mappings $[s]:A\to B$.

Let us define the category $C_\Theta(H)$.  Its objects have the
form $W/T$, where $W\in 0b\ \Theta^0$ and $T$ is an $H$-closed
congruence in $W$.  Morphisms of $C_\Theta(H)$ are the
homomorphisms of algebras.

It is proved that if $Var (H) = \Theta$ then the transitions $(X,
A) \to W(X)/A'$ and $W/T \to (X, T'_H)$ determine duality of the
categories $K_\Theta(H)$ and $C_\Theta(H)$.  In this case the
category $\Theta^0$ is a subcategory in $C_\Theta(H)$. The
skeleton of the category $K_\Theta(H)$ is denoted by $\widetilde
K_\Theta(H)$. This category is the category of {\it algebraic
varieties} over $H$. Correspondingly, the category $\widetilde
C_\Theta(H)$ is defined.

The category $K^0_\Theta (H) $ is always a subcategory in
$K_\Theta (H)$ \cite{Pl4}.

We consider also the categories $K_\Theta$ and $C_\Theta$ where
the algebra $H$ is not fixed. Correspondingly, we have the
categories $\tilde K_\Theta$ and $\tilde C_\Theta$.

{\bf 1.3} Consider a functor $Cl_H:\Theta^0 \to$ {\bf poSet},
where {\bf poSet} denotes the category of partially ordered sets.

This functor corresponds to every algebra $H$ in $\Theta$.  By
definition, for every $W \in 0b\ \Theta^0$ the poset $Cl_H(W)$ is
the set of all $H$-closed congruences $T$ in $W$ with the natural
ordering. Correspondingly, there is a lattice $Cl_H(W)$.

Let now a morphism
\[
s:W(Y) \to W(X)
\]
be given in $\Theta^0$.  It corresponds a map
\[
Cl_H(s): Cl_H(W(X)) \to Cl_H(W(Y)),
\]
defined by the rule $Cl_H(s)(T) = s^{-1}T$.
% corresponds to this $s$.
 Here $T \in Cl_H(W(X)); s^{-1} T$ is a congruence in $W(Y)$,
defined by the rule $w(s^{-1}T) w'$ if and only if $w^sTw'^s,\; w,
w' \in W(Y)$.  The congruence $s^{-1} T$ is also $H$-closed and
the mapping $Cl_H(s)$ is a morphism in the category  {\bf poSet}.

This defines a contravariant  functor $Cl_H$, which plays an
important part in the sequel.

%If $\Theta_1$ is a subvariety in $\Theta$, containing the algebra
%$H$, then there is also $Cl_H: \Theta^0_1 \to Set$, see [20].

In the same way one can consider a covariant functor
$Als_H:\Theta^0\to {\bf poSet}$, where $Als_H(W)$ is the poset of
algebraic sets in the affine space $\Hom(W,H)$.

\section{General look at the theory}

The main concepts of the theory are as follows:
\medskip

1. {\it Geometric properties of algebras $H$ in $\Theta$. An
algebra $H$ is considered in respect to its geometry and equations
over $H$}.
\medskip

2. {\it Geometric relations between algebras in $\Theta$}.
\medskip

3. {\it Structure of algebraic sets for every given algebra $H$
and every $W$. The lattice of algebraic sets in the given affine
space}.
\medskip

We will focus our attention on the problems related to Items 1 and
2. The item 3 is a separate topic which  requires the additional
clarity. %and information to be discussed from the point of view of
%the problems of generic type.

We now quote some working  notions   around which the theory
rotates.

First of all these are geometrical invariants of an algebra $H$:
{\it special categories and functors}. Categories are presented by
the {\it categories of algebraic sets and algebraic varieties}
$K_\Theta(H)$ and $\tilde K_\Theta(H)$. They are related to the
{\it categories $C_\Theta(H)$ and $\tilde C_\Theta(H)$}. Another
invariant of algebras is a contravariant functor $Cl_H : Var (H)^0
\to poSet$. Categories $K_\Theta$ and $C_\Theta$ are invariants of
the whole variety  $\Theta$.

The main features of algebras $H_i$ we are dealing with are {\it
geometrical noetherianity, logical noetherianity and geometrical
distributivity}. Relations between algebras are presented by the
notions of {\it geometrical equivalence, geometrical similarity,
geometrical compatibility, coincidence of geometries and
coincidence of lattices}. Here, coincidence of lattices in the
most general case is defined as isomorphism of functors of the
$Cl_{H_1} \to Cl_{H_2}\varphi$ type, where $\varphi$ is an
isomorphism of categories $\varphi : Var(H_1)^0 \to Var(H_2)^0$.

 Further we give all necessary  definitions.

\subsection{Examples of geometrical properties and relations}

{\bf 2.1.} {\it Geometrical  equivalence.}

\begin{defn}
Algebras $H_1$ and $H_2$ from
%in the given
$\Theta$ are called {\it geometrically equivalent} if for every
$W=W(X) \in Ob\; \Theta^0$ and every $T$ in $W$, we have
\[
T^{''}_{H_1} = T^{''}_{H_2}.
\]
This means also that $Cl_{H_1} = Cl_{H_2}$.
\end{defn}

   It is clear that if
the algebras $H_1$ and $H_2$ are geometrically equivalent, then
the categories $C_\Theta(H_1)$ and $C_\Theta(H_2)$ coincide.
Correspondingly,  the categories $K_\Theta(H_1)$ and
$K_\Theta(H_2)$ are isomorphic.

\begin{thm}
Algebras $H_1$ and $H_2$ are geometrically equivalent if and only
if
\[ LSC(H_1) = LSC(H_2).
\]
\end{thm}

Here the operator $L$ on classes of algebras is defined in the
usual local sense, i.e., for every class $\frak X$ an algebra $G$
belongs to $\frak X$, if every finitely generated subalgebra $H$
of $G$ belongs to $\frak X$. It can be proved that $LSC(\mathfrak
X)=\widetilde q Var (\mathfrak X)$, where $\widetilde q Var
(\mathfrak X)$ is the class of algebras which is determined by
infinitary quasiidentities of the class $\mathfrak X$.
Correspondingly, $ q Var (\mathfrak X)$ is the quasivariety which
is generated by the class $\mathfrak X$.
 Hence, geometrical equivalence of algebras means also that
\[
\tilde q Var (H_1) = \tilde q Var (H_2),
\]
i.e., $H_1$ and $H_2$ {\it have the same infinitary
quasiidentities}.
%\end{thm}

\medskip

%%%%%%%%%%%%%%%%%%%%%%%%%%%%%%%%%%%%%%%%%%%%%%%%%
{\bf 2.2.} {\it Geometrically and logically noetherian algebras}.
%%%%%%%%%%%%%%%%%%%%%%%%%%%%%%%%%%%%%%%%%%%%%%%%%%
\begin{defn}
 An  algebra $H\in\Theta$ is called {\it geometrically noetherian} if for an
arbitrary $W$ and $T$ in $W$ there exists a finite $T_0\subset T$
such that
 $$T^{''}_H = (T_0)^{''}_H.$$
\end{defn}
\medskip

%\begin{defn}
An algebra $H$ is  geometrically noetherian if and only if for
every $W$ and $T$ in $W$ there exists a finite subset $T_0\subset
T$ such that
%%%%
\[
\Big(\bigwedge\limits_{((w, w')\in T} (w\equiv w')\Big)
\Rightarrow w_0 \equiv w'_0
\]
holds in $H$ if and only if the quasiidentity
\[
\Big(\bigwedge\limits_{((w, w')\in T_0} (w\equiv w')\Big)
\Rightarrow w_0 \equiv w'_0
\]
holds in $H$. Here $T_0$ is independent from $(w_0,w'_0)$.
%\end{defn}

\begin{defn}
In case when $T_0$ depends on $(w_0,w'_0)$ we call $H$ {\it
logically noetherian}.
\end{defn}

 The notion to be logically noetherian means
also that $T''$ coincides with $\bigcup T_\alpha''$ where the
union is taken over all finite sunsets $T_\alpha$ in $T$.

Obviously, if $H$ is geometrically noetherian, then $H$ is
logically noetherian.

An algebra $H$ turns to be  {\it geometrically noetherian} if and
only if in every $W=W(X)$ the ascending chain condition for
$H$-closed congruences holds. Dually, the descending chain
condition for algebraic sets in $\Hom(W(X), H)$ holds in
geometrically noetherian algebras. An algebra $H$ is logically
noetherian if the union of  a directed set of $H$-closed
congruences is also $H$-closed.

\begin{thm} \cite{MR} Let $H_1$ and $H_2$ be logically noetherian algebras.
 They are
geometrically equivalent if and only if $q Var (H_1) = q Var
(H_2)$.
\end{thm}

%It follows from  Theorem 2 that the equality $LSC(H) = qVar (H)$
%holds if and only if $H$ is logically noetherian.  This together
%with the presentation for $qVar (H)$ imply

\begin{thm} \cite{MR}
If $H\in\Theta$ is not logically noetherian, then there exists an
ultrapower $H'$ of $H$ such that the algebras $H$ and $H'$ are not
geometrically equivalent.

 However, these algebras have the same elementary theories and, in
particular, the same quasiidentities.
\end{thm}

These theorems lead to the following general problem:
\medskip

{\it For which varieties $\Theta$ there exist non-logically
noetherian algebras in $\Theta$?

How often these algebras can appear?}
\medskip

The existence of such phenomenon for groups is proved in the paper
by K.Gobel, S.Shelah \cite{GS}. The idea of their proof is based
on the existence of the continuum different 2-generated simple
groups \cite{LS}. For representations of groups the result is
proved by A.Tsurkov \cite{PT}. For associative algebras over a
field the result also holds \cite{Pl5}. In the recent paper of
Lichtman - Passman  \cite{LP} the existence of the continuum of
3-generated simple algebras is proved.

The results and notions above are of universal character. In
particular, they can be applied to multi-sorted algebras. Further
we consider concrete $\Theta$ and mostly for them we formulate
problems.

\section{Geometrical properties of algebras. Problems}
\medskip

{\bf Problem 1}. {\it Let $G=A wr B$ be a wreath product of some
groups $A$ and $B$.
\medskip

1. When $G$ is geometrically noetherian?
\medskip

2. When $G$ is  logically noetherian but not geometrically
noetherian?
\medskip

3. Are there groups $G=AwrB$ which are not logically noetherian
for some  appropriate $A$ and $B$. }

\medskip

It is known that any free group $W(X)$ is geometrically noetherian
(Guba \cite{G}). Moreover, every group or algebra which admits
faithful finite dimensional representation is geometrically
noetherian (Miasnikov-Remeslennokov \cite{MR}, Kanel-Belov). Every
finite dimensional representation of a group is geometrically
noetherian (Tsurkov \cite{PT},

{\bf Problem 2}.
\medskip

{\it Is it true that every free Lie algebra $W(X)$ is
geometrically noetherian? }
\medskip

 Most likely the answer is negative.
Thus arises the following:
\medskip

{\bf Problem 3}.
\medskip

{\it Is it true that every free Lie algebra $W(X)$ is logically
noetherian?}

\medskip

{\bf Problem 4}.
\medskip

{\it Is it true that every free associative algebra $W(X)$ is
geometrically noetherian?}

Here the expected answer is also seems to be negative. Then:

\medskip

{\bf Problem 5}.
\medskip

{\it Is it true that every free associative algebra $W(X)$ is
logically noetherian?}

\medskip

Any two free groups have the same quasi-identities. The similar
fact is valid for free associative and free Lie algebras. Free
groups are also geometrically noetherian. Geometrical
noetherianity together with coincidence of quasi-identities
implies that any two free groups are geometrically equivalent.
Thus, the free groups have the same logic of quasi-identities and
the same geometry. The positive solution of Problems 3 and 5 would
mean that the same fact holds true for the free Lie algebras and
free associative algebras.

{\bf Problem 6}.
\medskip

{\it Is it true that there exists a continuum of different
$k$-generated simple Lie algebras? Here $k$ is fixed.}

\medskip

{\bf Problem 7}.
\medskip

{\it Is it true that there exists a non-logically noetherian Lie
algebra?}

\medskip

The next problems are devoted to lattices of algebraic sets.

\begin{defn}
An algebra $H$ is called {\it geometrically distributive} if for
every $W$ the lattice of algebraic sets $Als_H(W)$ (and the
lattice $Cl_H(W)$, respectively) is distributive.
\end{defn}

The {\it geometrically modular} algebras are defined in the
similar way.

\medskip

{\bf Problem 8}.
\medskip

{\it Which algebras $H$ are geometrically distributive?}
\medskip

This problem makes sense for groups, groups with the fixed group
of constants, and for other varieties $\Theta$.

\medskip

{\bf Problem 9}.
\medskip

{\it Which algebras $H$ are geometrically modular?}
\medskip

We introduced earlier the category $K_\Theta$ of algebraic sets
without the fixed set $H$.

\medskip

{\bf Problem 10}.
\medskip

{\it When the categories ${K_{\Theta_1}}$ and ${K_{\Theta_2}}$ are
isomorphic and when they are equivalent? Consider separately the
case when $\Theta_1$ and $\Theta_2$ are subvarieties of a bigger
variety $\Theta$.}

This problem should be related with the known results of McKenzie
\cite{Mc} about equivalence of two varieties of algebras.

\section{Other geometrical relations between algebras}

We have defined the notion of geometric equivalence of algebras
$H_1$ and $H_2$. Now we define two more general notions.

{\bf 4.1.} First we recall the definition of isomorphism of
functors.

Let two functors $\varphi_1,\varphi_2: C_1\to C_2$ of the
categories $C_1,C_2$ be given. A {\it homomorphism (natural
transformation) of functors} $s:\varphi_1\to \varphi_2$ is a
function, relating a morphism in $C_2,$ denoted by $s_A:
\varphi_1(A)\to \varphi_2(A)$ to every object $A$ of the category
$C_1.$ For every $\nu: A\to B$ in $C_1$ there is a commutative
diagram
  % \begin{CD}
%$$\CD \Hom(F(X),H) @>\tilde s>> \Hom(F(X^0),H) \\ @V\alpha_X VV @VV\alpha_{X^0} V\\
%H^n@<s^\alpha<< H^k\endCD$$ \noindent
%%%%%%%%%%%%%%%%%%%%%%%%%%%%%%%%%%%%%%%%%%%%%%%%%%%%%%%%%%%%%%

$$
\CD
\varphi_1(A) @> s_A>> \varphi_2(A)\\
@V\varphi_1(\nu)VV @VV\varphi_2(\nu) V\\
\varphi_1(B) @>s_B>> \varphi_2(B)
%\end{CD}
\endCD
$$
\medskip

\noindent
in the case of covariant $\varphi_1$ and $\varphi_2.$
For contravariant $\vp_1$ and $\vp_2$ the corresponding diagram is
$$\CD
%\begin{CD}
\varphi_1(B) @> s_B>> \varphi_2(B)\\
@V\varphi_1(\nu) VV @VV\varphi_2(\nu) V\\
\varphi_1(A) @>s_A>> \vp_2(A)
%\end{CD}
\endCD
$$
\medskip

\noindent
 An invertible  homomorphism $s:\vp_1\to\vp_2$ is called {\it an isomorphism
(natural isomorphism) of functors}. The isomorphism property holds
if $s_A:\vp_1(A)\to\vp_2(A)$ is an isomorphism in $C_2$ for any
$A$.

\begin{defn}
Algebras $H_1$ and $H_2$ are called {\it geometrically similar},
if there exists an automorphism $\varphi:\Theta^0\to\Theta^0$ such
that there is a {\it correct isomorphism} of functors
$\alpha(\varphi): Cl_{H_1}\to Cl_{H_2}\varphi$.
\end{defn}

 Here, correctness
means compatibility with the automorphism $\varphi$. Namely, let
$s_1,s_2: W_1\to W_2$ be given, and $T$ be $H_1$-closed in $W_2$.
Denote, $T^*=\alpha(\varphi)_{W_1}(T)$. There are canonical
homomorphisms $\mu_T:W_2\to W_2/T$ and $\mu_{T^*}:\varphi(W_2)\to
\varphi(W_2)/T^*$. Correctness means that $\mu_Ts_1=\mu_Ts_2$
holds if and only if
$\mu_{T^*}\varphi(s_1)=\mu_{T^*}\varphi(s_2).$

\begin{defn}
Algebras $H_1$ and $H_2$ are called {\it geometrically
compatible}, if there exists an autoequivalence of the category
$\Theta^0\ _\psi\leftrightarrows^\varphi\Theta^0$ such that there
are the natural transformations of functors
$$
\alpha(\varphi): Cl_{H_1}\to Cl_{H_2}\varphi,
$$
$$
\alpha(\psi): Cl_{H_2}\to Cl_{H_1}\psi,
$$
which are compatible as before with $\varphi$ and $\psi$.
\end{defn}

%\begin{defn}
%Algebras $H_1$ and $H_2$ are called {\it geometrically
%compatible}, if for some autoequivalence of the category
%$\Theta^0\ _\psi\leftrightarrows^\varphi\Theta^0$ we have the
%natural transformations of functors
%$$
%\alpha(\varphi): Cl_{H_1}\to Cl_{H_2}\varphi,
%$$
%$$
%\alpha(\psi): Cl_{H_2}\to Cl_{H_1}\psi.
%$$
%which are compatible as before with $\varphi$ and $\psi$.
%\end{defn}

{\bf 4.2.} We consider also {\it the correct isomorphisms}  of the
categories $C_\Theta(H_1) \to C_\Theta(H_2)$. These are
isomorphisms which induce an automorphism of the category
$\Theta^0$. The correct isomorphisms of the categories
$K_\Theta(H_1) \to K_\Theta(H_2)$ are defined in a similar way.

\begin{thm} Suppose $\Var H_1=\Var H_2=\Theta$. The categories $K_\Theta(H_1)$ and $K_\Theta(H_2)$
are correctly isomorphic if and only if $H_1$ and $H_2$ are
geometrically similar.
\end{thm}

\begin{thm} Suppose $\Var H_1=\Var H_2=\Theta$. The categories $K_\Theta(H_1)$ and $K_\Theta(H_2)$
are correctly equivalent if and only if $H_1$ and $H_2$ are
geometrically compatible.
\end{thm}

Correctness here means also that the lattices of algebraic sets
for $H_1$ and $H_2$ are the same.

These theorems are of the universal character. Each of them should
be specified for particular varieties $\Theta$. This
specialization depends very much on the description of the group
$\Aut(\Theta^0)$.
\medskip

{\bf Problem 11}
\medskip

{\it Consider the similar problems  without assumption of
correctness for isomorphisms and equivalences of categories.}

\section{$\Aut(\Theta^0)$}

All automorphisms of the category $\Theta^0$ are known in the
following cases (see \cite{Be},\cite{MPP1},\cite{MPP2},
\cite{Msh}, \cite{LiP}, \cite{Pl5}, \cite{MS})
.

\medskip

{\it

 1. Groups.
\medskip

2. Groups with a free group of constants, $\Grp-F$.
\medskip

3. Associative and commutative algebras, $\Com-P$.
\medskip

4. Associative algebras.
\medskip

5. Lie algebras.
\medskip

6. $K$-modules, $K$ is an   $IBN-$ ring.
\medskip

7. Semigroups.
\medskip
}

 In the situation of Lie algebras the description of
automorphisms uses the description of the group
$\Aut(\End(W(x,y))$. This observation motivates the following:
\medskip

{\bf Problem 12}
\medskip

 {\it Study the
group $\Aut(\End(W(X))$, where $W(X)$ is the free Lie algebra over
a finite set $X$.}

\medskip

{\bf Problem 13}
\medskip

{\it Study the group $\Aut(\Theta^0)$ for various interesting
subvarieties of the variety of all groups. For example for the
varieties $\frak N_c$, $ \frak A^2$, etc. }

\medskip

{\bf Problem 14}
\medskip

{\it Study the group $\Aut(\Theta^0)$ for various interesting
subvarieties of the variety of all Lie algebras. }

\medskip

{\bf Problem 15}
\medskip

{\it Study the group $\Aut(\Theta^0)$ for various interesting
subvarieties of the variety of all associative algebras. }

\section{Algebras with the same algebraic geometry}

Recall that we look at the notion of coincidence of geometries in
two variants.
\medskip

1. The categories $K_\Theta(H_1)$ and $K_\Theta(H_2)$ are
isomorphic.
\medskip

2. The categories $K_\Theta(H_1)$ and $K_\Theta(H_2)$ are
equivalent.
\medskip

In fact, the second case means that the categories of algebraic
varieties  $\widetilde K_\Theta(H_1)$ and $\widetilde
K_\Theta(H_2)$ are isomorphic.

We consider the specific varieties $\Theta$. The problem of
coincidence of geometries is solved in the following cases:

\medskip

{\it
 1. For the classical algebraic geometry \cite{Be}.
\medskip

2. For the non-commutative algebraic geometry related to the

variety of all associative algebras \cite{MPP2}.
\medskip

3. For the algebraic geometry in the variety of all Lie algebras
\cite{Pl5}.
\medskip

4. For the geometry in the variety of all groups \cite{Pl5}.
\medskip

5. For the variety $\Grp-F$ \cite{BPP}.
\medskip
}

 {\bf Problem 16}.
\medskip

 {\it Investigate coincidence of geometries for some subvarieties
of the variety of all groups. }

\medskip

{\bf Problem 17}.
\medskip

{\it Investigate coincidence of geometries for some subvarieties
of the variety of all  Lie algebras. }

\medskip

{\bf Problem 18}.
\medskip

{\it Investigate coincidence of geometries for some subvarieties
of the variety of all associative algebras. For example, for the
subvariety $\Theta$  given by a single polynomial identity. }

\medskip

Solution of these problems heavily depends on the solution of the
problem of the group $\Aut(\Theta^0)$ description.

\section{Coincidence of the lattices of algebraic sets}

We consider the following variants for the definition of lattices
coincidence.

\medskip

1. Coincidence of the functors $\Cl_{H_1}$ and $\Cl_{H_2}$.

\medskip

2. Isomorphism of the functors $\Cl_{H_1}$ and $\Cl_{H_2}$.

\medskip

3. The functor $\Cl_{H_1}$ is isomorphic to $\Cl_{H_2}\varphi$,
where $\varphi$ is an automorphism of the category $\Theta^0$.

\medskip

In the first case the algebras $H_1$ and $H_2$ are geometrically
equivalent and the lattices in $W$ corresponding to $H_1$ and
$H_2$ coincide.

In the second case an isomorphism of functors $\alpha:\Cl_{H_1}\to
\Cl_{H_2}$ provides an isomorphism of the corresponding lattices
for every $W$. Besides, there is a compatibility with the
morphisms.

In the third case for every $W$ there is an isomorphism of the
lattices $\Cl_{H_1}(W)$ and $\Cl_{H_2}(\varphi(W))$.
\medskip

{\bf Problem 19}
\medskip
{\it For which algebras $H_1$ and $H_2$ there is an isomorphism of
the functors $\Cl_{H_1}$ and $\Cl_{H_2}.$ }

\medskip

{\bf Problem 20}
\medskip
{\it  For which algebras $H_1$ and $H_2$ there is an isomorphism
between $Cl_{H_1}$ and $\Cl_{H_2}\varphi.$} for some
$\varphi:\Theta^0\to\Theta^0$.

\medskip

If the algebras $H_1$ and $H_2$ are geometrically similar then
such an isomorphism exists. Thus, if $H_1$ and $H_2$ have the same
geometries then the corresponding lattices coincide. The converse
statement is not true and this makes everything more attractive.
%%%%%%%%%%%%%%%%%%%%%%%%%%%%%%%%%%%%%%%%%%%%%%%%%%%%%%%%%%%%%%%%%%%%%%%%

The problems above seem new also for the classical situation
$Com-P$, where $L_1$ and $L_2$ are two extensions of a ground
filed $P$.

In particular, what can be said about $L_1$ and $L_2$ if for every
$W$ the lattices $\Cl_{L_1}(W)$ and $\Cl_{L_2}(W)$ are isomorphic?

Coincidence of these lattices means that $L_1$ and $L_2$ are
geometrically equivalent. Hence, in this case the logics of
quasiidentities for $L_1$ and $L_2$ are the same. But we are
interested in conditions providing isomorphism of lattices.

%%%%%%%%%%%%%%%%%%%%%%%%%%%%%%%%%%%%%%%%%%%%%%%%%%%%%%%%%%%%%%%%%%%%%%%%

\section{Representations}

{\bf 8.1.} Let $K$ be a commutative, associative ring with unit.
We consider the category-variety $\Theta=\Rep-K$. General
references on the theory in question are
\cite{PV},\cite{Pl3},\cite{V1}, \cite{V2}.

Objects of this category are representations $(V,G)$, where $V$ is
a $K$-module and $ G$ is a group acting on $V$. These $(V,G)$ are
two-sorted algebras.

The action $G$ on $V$ is denoted by $\circ$ and for every $a\in V$
and $g\in G$ we have $a\circ g\in V$. The action $\circ$ satisfies
the natural identities.

Morphisms in $\Theta=\Rep-K$ have the form
$$
\mu=(\alpha,\beta): (V_1,G_1)\to (V_2,G_2)
$$
where $\alpha\in \Hom_K(V_1,V_2),$ $\beta \in \Hom_(G_1,G_2),$ and
$(a\circ g)^\alpha=a^\alpha\circ g^\beta.$

$\Ker \mu=(\Ker \alpha, \Ker\beta)=(V_0,H)$ is a congruence in
$(V_1,G_1)$ in the following sense: $H_0$ is $G_1-$invariant
submodule in $V_1$ and $H$ acts trivially in $V_1/V_0$. We have
the factor-representation $(V_1,G_1)/(V_0,H)=(V_1/V_0,G/H)$ with
the natural theorem on homomorphisms. For a given set
$\mu_i=(\alpha_i,\beta_i):(V_1,G_1)\to (V_2,G_2)$ we have
$\bigcap\Ker\mu_i=(\bigcap\Ker\alpha_i, \bigcap\Ker\beta_i)$.

Free objects $W$ in the category $\Theta$ are denoted by
$W=W(X,Y)$, where $X$ and $Y$ is a pair of sets. Here:
$$
W(X,Y)=(XKF(Y),F(Y))
$$
where $F=F(Y)$ is the free group over $Y$, $KF$ is the group
algebra, $XKF$ is the free $KF$-module over the set $X$.

For every $w\in XKF$, $w=x_1u_1+\ldots+x_nu_n$, $u_i\in KF,$ and
$f\in F$ we have
$$
w\circ f=x_1(u_1f)+\ldots+x_n(u_nf).
$$

This is the free representation in the categorical sense over the
two-sorted set $(X,Y)$. The two-sorted equality $w\equiv 0$ is
considered as an {\it action-type equality}, while $f\equiv 1$ is
a group equality.

In the book \cite{PV} the action-type varieties of representations
have been considered. These varieties lie in $\Rep-K$ and can be
defined by identities of the type $x\circ u$.

{\bf 8.2.} In $\Rep-K$ one  can consider the different operations:
Cartesian--direct products, free products--coproducts,
subrepresentations,  quotients, etc. The following two
constructions are of the special type.

{\it Triangular products.} For the given representations
$(V_1,G_1)$ and $(V_2,G_2)$ consider their triangular product
$(V_1,G_1)\triangledown(V_2,G_2).$ This is the representation
$(V_1+V_2,G)$, where $g\in G$ has the form $$
     g=\left[
\begin {array}{ccc}
g_2&\varphi g_1\\
\noalign{\medskip} 0& g_1
\end {array}
\right]= \left[
\begin {array}{ccc}
1&\varphi \\
\noalign{\medskip} 0& 1
\end {array}
\right] \left[
\begin {array}{ccc}
g_2& 0\\
\noalign{\medskip} 0& g_1
\end {array}
\right],
 $$
$g_1\in G_1$, $g_2\in G_2$, $\varphi \in \Hom(V_2,V_1).$

For $a\in V_1, $ $b\in V_2$ we have $ a\circ g=a\circ g_1$; $
b\circ g= b\circ g_2+(b\varphi)\circ g_1$. Here, $(V_1,G)$ is
related to $(V_1,G_1)$ and $(V_1+V_2/V_1,G)$ to $(V_2,G_2).$

We consider also wreath product $(V,H)wr G= (V^G,Hwr G)$.

{\bf 8.3.} In $\Rep-K$ we consider varieties of the general form
and action-type varieties. For the last ones the semigroup
$\mathfrak M$ of such the varieties $\mathfrak X$ is treated. The
multiplication in $\mathfrak M$ is defined by the rule: $(V,G)\in
\mathfrak X_1 \mathfrak X_2$ if for some $G$-invariant submodule
$V_0\subset V$ we have $(V_0,G)\in \mathfrak X_1, $ $(V/V_0,G)\in
\mathfrak X_2$.

If $K$ is a field then the semigroup $\mathfrak M$ is a free
semigroup and we have here
$$
\Var((V_1,G_1)\triangledown(V_2,G_2))=\Var(V_1,G_1)\Var(V_2,G_2).
$$

Let now $\mathfrak N$ be the semigroup of group varieties. For
$\mathfrak X\in \mathfrak M$ and $\Theta \in \mathfrak N$ we
consider the product $\mathfrak X \times \Theta\in\mathfrak M$
defined as follows:  $(V,G)\in \mathfrak X \times \Theta$ if for
some invariant subgroup $H$ in $G$ we have $(V,H)\in \mathfrak X$,
$G/H\in\Theta$. Now $\mathfrak N$ acts in $\mathfrak M$ as a
semigroup of endomorphisms of $\mathfrak M.$

The principal theorem in this theory says that the action of
$\mathfrak N$ in $\mathfrak M$ is free. This means also that every
$\mathfrak X$ in $\mathfrak M$ can be uniquely presented in the form
$$
\mathfrak X=(\mathfrak X_1\times\Theta_1)\ldots(\mathfrak
X_n\times\Theta_n),
$$
where all $\mathfrak X_i$ are irreducible.

The triangular products and wreath products play the crucial role
in the proof of the theorem above.

\section{Algebraic geometry in representations}

{\bf 9.1.} We consider $\Hom(W,(V,G))$ as the {\it affine space}
over the given representation $(V,G).$  Here, $W=W(X,Y)$ is the
free representation over the finite sets $X$ and $Y$. {\it Points}
here are homomorphisms
$$
\mu : W\to (V,G).
$$
Take $T=(T_1,T_2)$, where $T_1$ is a set of action-type equalities
in $W$ and $T_2$ is a set of group equalities.

Denote

\[
\left\{%\begin{matrix}
\begin{aligned}&T'_{(V,G)} =A=\{ \mu=(\alpha,\beta):W\to (V,G)\ | T_1\subset \Ker\alpha\ ,
\ T_2\subset \Ker\beta\}\\
&A'_W =T=\bigcap\limits_{\mu\in A} \Ker
\mu=(\bigcap\limits_{\alpha}\Ker\alpha,
\bigcap\limits_{\beta}\Ker\beta)
%\end{matrix}
\end{aligned} \right.
\]
$$
T_1=\bigcap\Ker\alpha\ , \ T_2=\bigcap\Ker\beta.
$$

Let $\Id_G(F)$ be the verbal subgroup of all identities of the
group $G$ in $F=F(Y)$. In every case we have $\Id_G(F)\subset
T_2$.

A set  $A$ of the form $A=T'$ is {\it an algebraic set}, and
$T=A'$ is a $(V,G))$-closed congruence in $W$. The definitions
above specialize the
general definitions from universal algebraic
geometry for the case of the variety-category of group
representations. Some results of the universal algebraic geometry
in multi-sorted $\Theta$ are applicable in this case.

{\bf 9.2.} Now we consider action-type AG in representations.

For the given $W=W(X,Y)=(XKF,F)$ $F=F(Y)$ we take a set $T\subset
XKF$. We view $T$ as a set of action-type equalities.

Denote

\[
\left\{%\begin{matrix}
\begin{aligned}&T^v =A=\{ \mu=(\alpha,\beta):W\to (V,G)\ | T\subset
\Ker\alpha\}\\
&A^v=T=\bigcap\limits_{\alpha}\Ker\alpha.
%\end{matrix}
\end{aligned} \right.
\]

Here $A=T^v$ is an {\it action-type algebraic set}, and $T=A^v$ is
an {\it an action-type $(V,G)$-closed $F$-invariant submodule in
$XKF$}.

It is easy to see that an algebraic set $A$ is an action-type
algebraic set if and only if all points of the type $(0,\beta)$
belongs to $A$. From this follows that if $A$ is an action-type
algebraic set, then
$$
A'=(A^v,\Id_G(F)).
$$

As before one can define  the notions of the geometrically
equivalent representations as well as the notions of geometrically
and logically noetherian representations. These definitions refer
to general case and also to action-type case.

We can consider also the categories $K_\Theta(V,G)$ and
$C_\Theta(V,G)$ for the general case and the categories
$K^{at}_\Theta(V,G)$ and $C^{at}_\Theta(V,G)$ for the action-type
case. Here, $\Theta=\Rep-K$ or a subvariety of $\Theta=\Rep-K$ .

{\bf 10.3.} Once again the open problems.

\medskip

 {\bf Problem 21}
\medskip
{\it  When the representations $(V_1,G_1)$ and $(V_2,G_2)$ have
the same geometry}.

\medskip

This question relates to general situation and to action-type
case.

With this problem the notions of geometrically similar and
geometrically compatible representations are associated.
Furthermore, the latter notions are connected with the
automorphisms and autoequivalences of the category
$\Theta^0=(\Rep-K)^0$.

So, we have the following

\medskip

 {\bf Problem 22}
\medskip

{\it  Investigate the group $\Aut (\Rep-K)^0$.}

\medskip

Recall that an automorphism of the category $C$ is called {\it
inner} if it is isomorphic to unity automorphism $1_C$.

Inner automorphisms constitute an invariant subgroup in the group
$\Aut$$ (\Rep-K)^0$.

%%%%%%%%%%%%%%%%%%%%%%%%%%%%%%%%%%%%%%%%%%%%%%%%%%%%%%%%%%%%%%%%%%%%%%
One can speak also on semi-inner automorphisms. In their definition
the automorphisms $\sigma$ of the ring $K$ take part. These
automorphisms  form a subgroup in the group $\Aut (\Rep-K)^0$. There
is also a special mirror (anti-inner) automorphism $\delta$. The
existence of such automorphism is based on the consideration of the
opposite group and opposite representation.

The problem is to prove that the group  $\Aut (\Rep-K)^0$ is
generated by automorphisms above.

%%%%%%%%%%%%%%%%%%%%%%%%%%%%%%%%%%%%%%%%%%%%%%%%%

Repeating the arguments from \cite{Pl5} one can prove that if the
similarity of two representations $(V_1,G_1)$ and $(V_2,G_2)$  is
related to an inner automorphism $\varphi$ of the category
$(\Rep-K)^0$, then the representations are geometrically
equivalent in general. This implies that they are also action-type
equivalent.

\medskip

 {\bf Problem 23}
\medskip

{\it  What is the relation between the representations $(V_1,G_1)$
and $(V_2,G_2)$  if they are geometrically similar and the
similarity is based on an semi-inner automorphism $\varphi$ of the
category $(\Rep-K)^0$?}

\medskip

This problem is connected with

\medskip

 {\bf Problem 24}
\medskip

{\it  Investigate the group $\Aut (\End(KF,F))$.}

\medskip

Let us discuss this problem in more detail. For every
representation $(V,G)$ we have the group $\Aut (\End(V,G))$. Let
$\xi=(s,\tau)$ be an invertible element of the semigroup
$\End(V,G)$. Then $\xi$ is also an automorphism of the
representation $(V,G)$. An inner automorphism $\hat\xi$ of the
semigroup $\End(V,G)$ corresponds to it.

For every $\mu=(\alpha,\beta)\in \End(V,G)$ we have
$$
\hat\xi(\mu)=\xi\mu\xi^{-1}=(s\alpha s^{-1},\tau\beta\tau^{-1}).
$$
All these $\hat\xi$ form a normal subgroup in the group $\Aut
(\End(V,G))$.

Consider the pairs $\xi=(s,\tau)$, where $s$ is a
semi-automorphism of the $K$-module $V$ and $\tau\in\Aut(G)$.
There is $\sigma\in \Aut(K)$ such that for every $\lambda\in K$
and $a\in V$ we have:
$$ s(\lambda a)=\lambda^\tau s(a).
$$
Besides, $(a\circ g)^s=a^s\circ g^\tau$ for every $a$ and $g$.
Here, $\xi$ is a semi-automorphism of the representation $(V,G)$,
it does not belong to the semigroup $\End(V,G)$. However, $\xi$
induces automorphism $\hat \xi$ of this semigroup. Here, $\hat
\xi$ is a semi-inner automorphism of the semigroup $\End(V,G)$ and
all these automorphisms form a subgroup in $\Aut(V,G)$.

Let, further, $\varphi$ be an arbitrary automorphism of the
semigroup $\End(V,G)$. For every $\mu=(\alpha,\beta)\in \End(V,G)$
we have
$$
\varphi(\mu)=(\varphi_1(\mu),\varphi_2(\mu)).
$$
Here, $\varphi_1:\End(V,G)\to \End V$, $\varphi_2:\End(V,G)\to
\End V$ are the homomorphisms of semigroups.

One can calculate the conditions on homomorphisms $\varphi_1$ and
$\varphi_2$ which provide the homomorphism $\varphi$ of the
semigroup $\End(V,G)$. However, it is not clear how to deduce from
these conditions the "real constructions".

All above can be applied to the representations of the kind $(KG,G)$
and, in particular, to $(KF,F)$. In this important case one has to
take into account the theorem of Formanek \cite{For}, which says
that all automorphisms of the semigroup $\End(F)$ are inner.  Is it
possible to state that all automorphisms of the semigroup
$\End(KF,F)$ are semi-inner or of the form  $\varphi \delta $ where
$\varphi$ is semiinner.? Or one can construct a counter example?

Now consider the problems of the different kind. All these
problems should be stated separately for the general and
action-type cases.

\medskip

 {\bf Problem 25}
\medskip

{\it  Consider the representations $(V_1,G_1)$ and $(V_2,G_2)$
from the point of view of coincidence of the corresponding
lattices of algebraic varieties.}

\medskip

\medskip

 {\bf Problem 26}
\medskip

{\it Is it true that the  representation $(XKF,F)$ is
geometrically noetherian or logically noetherian?}

\medskip

Recall here that the group $F$ is geometrically noetherian.

\medskip

 {\bf Problem 27}
\medskip

{\it Consider the notions of geometrical and logical noetherianity
in respect to triangular products of representations and wreath
products  of a representation and a group.}

\medskip

{\bf Problem 28}
\medskip

Let $G$ be not logically noetherian group. Whether the group algebra
$PG$ is also not logically noetherian for some field $P$. The same
question is meaningful for the regular representation $(PG,G)$ in
action-type geometry.

\end{document}